\documentclass[12pt,a4paper]{article}
\usepackage[]{fontenc}
\usepackage[latin1]{inputenc}
 \usepackage[colorlinks]{hyperref}
\usepackage{a4,amsmath,amssymb,amsfonts,amsthm}
\usepackage{amsrefs}
\usepackage[all]{xy}
\newcommand{\legendre}[2]{\genfrac{(}{)}{}{}{#1}{#2}}

\def\mR{\mathbb{R}}
\def\mZ{\mathbb{Z}}
\def\mQ{\mathbb{Q}}
\def\mC{\mathbb{C}}

\def\Fp{{\mathfrak p}}

\def\Spec{{\rm Spec}}

 \def\Gal{{\rm Gal}}
 
 \def\det{{{\rm det}}}

  \def\Spec{{\rm Spec}}
\def\Id{{\rm Id}}

\def\Cl{{\rm Cl}}

\def\beginProof{\par{\bf Proof. }}
 \def\endProof{${\qed}$\par\smallskip}

 \def\to{\rightarrow}
 
 \def\mFa{{\mathfrak a}}
 
  \def\mFb{{\mathfrak b}}

 \def\mQ{{\Bbb Q}}
 \def\mR{{\Bbb R}}
 \def\mZ{{\Bbb Z}}
 \def\mC{{\Bbb C}}

 \def\CH{{\mathcal H}}
 \def\CO{{\mathcal O}}
 \def\CA{{\mathcal A}}
 
 \def\CB{{\mathcal B}}
  \def\CH{{\mathcal H}}

 \def\md{{\rm mod\ }}
 \def\refeq#1{(\ref{#1})}

\def\P1{{{\bf P}^1}}

\def\without{\backslash}

\def\Lie{{\rm Lie}}

\def\Pet{{\rm Pet}}

\def\Tau{{\mathcal T}}

\def\c1{{\rm c}_1}

 \newtheorem{theor}{Theorem}[section]
 \newtheorem{prop}[theor]{Proposition}
 
 \newtheorem{cor}[theor]{Corollary}
 \newtheorem{lemma}[theor]{Lemma}

 \newtheorem{conj}[theor]{Conjecture}

\newtheorem{rem}[theor]{Remark}

\begin{document}

 \author{Vincent Maillot\footnote{Institut de Math\'ematiques de Jussieu, Universit\'e Paris 7 Denis Diderot, C.N.R.S., Case Postale 7012, 2 place
Jussieu, F-75251 Paris Cedex 05, France, E-mail: vmaillot@math.jussieu.fr}\,\, and Damian R\"ossler\footnote{Mathematical Institute, University of Oxford, Andrew Wiles Building, Radcliffe Observatory Quarter, Woodstock Road, Oxford, OX2 6GG, United Kingdom, E-mail: damian.rossler@maths.ox.ac.uk}}

 \title{The conjecture of Colmez and reciprocity laws for modular forms}
\maketitle
\begin{abstract}
In \cite{Colmez-Periodes}, P. Colmez formulated a remarkable conjecture, which asserts that the Faltings height of a CM abelian variety can be computed as a linear combination of logarithmic derivatives of Artin $L$-functions. Noting that the Faltings height is an average of transcendental quantities summed over the embeddings of a number field of definition of the abelian variety, we propose a refinement of this conjecture, which identifies each of these transcendental quantities. We also show how our conjecture would imply the existence of fine reciprocity laws for Siegel modular forms with rational coefficients evaluated at CM points, and we prove our conjecture for elliptic curves, using old results of Siegel and Hasse on elliptic units.
\end{abstract}

 \tableofcontents

\section{Introduction}

In \cite{Colmez-Periodes}, Colmez established a remarkable explicit formula for the Faltings height 
of an abelian variety which has complex multiplications by the ring of integers of a CM field 
whose normalisation is an abelian extension of $\mQ$. This formula is 
a higher dimensional generalisation of a classical formula of Chowla and Selberg, 
which was first established by Lerch (see \cite{Lerch}, \cite{CS-Epstein}). 
Colmez also conjectured that a very natural generalisation of his formula is true without any assumption on the Galois group of the CM field.
There was a lot of work around his formula and its conjectural extension over the last thirty years. It was in particular shown in \cite{KR-IV} that Arakelov theory could be used to give a new proof of a weak form of his formula, and in more recent years, an averaged form of his conjecture was established (see \cite{Yuan-Zhang-AC}, \cite{AGHP-Fal}). 

The aim of the present text is to formulate a refinement of the conjecture of Colmez (Conjecture \ref{mainC} below), which is 
suggested by Arakelov theory. We will also explain how one can deduce our conjecture from the conjecture of Colmez in any situation where one can find a Siegel modular form which satisfies a reciprocity law and which has suitable integrality properties with respect to $A$ and $E$ (see Lemma \ref{lemAUTR}). Using classical results of analytic number theory going back to Hasse and Siegel, we show that in the case of CM elliptic curves, the discriminant modular form has the required properties, thus proving our conjecture for CM elliptic curves (see Proposition \ref{propELL}).  

Section \ref{secSTAT} below contains the statement of our conjecture. 
Section \ref{secSIEG} is about Siegel modular forms and what the non vanishing of such a form 
at a CM point means in the context of our conjecture.
Section \ref{secELL} establishes our conjecture for CM elliptic curves.

\section{The statement}
\label{secSTAT}

Before formulating our conjecture, we will recall Colmez's original conjecture.

Let $A$ be an abelian variety of dimension $d$ defined over a number field $K$. 
Suppose that the N\'eron model $\cal A$ of 
$A$ over 
${\cal O}_{K}$ has semi-stable reduction at all the places of $K$.
 Let $\omega_{\cal A}$ be the ${\cal O}_{K}$-module, which consists 
of $\CA$-invariant global sections of the sheaf $\Omega^d_{\CA/S}$ of relative $d$-forms on $\CA$. 
If $\sigma:K\to\mC$ is an embedding of ${K}$ in $\bf C$, we let ${A}_{\sigma}(\mC)$ be the manifold of complex points of the variety $A_\sigma:={A}\times_{\sigma({K})}\mC$. 
Let $\alpha\in\omega_\CA$. We write $\alpha_\sigma\in\Omega^{d}_{{A}_{\sigma}}$ for the  the differential form obtained by base change by $\sigma:K\to\mC$. Abusing language, 
we also write $\alpha_\sigma$ for the corresponding analytic differential form on ${A}_{\sigma}(\mC).$

The
{\it Faltings height} (or {\it modular height}) of $A$ is the quantity
\begin{equation}
h_{\rm Fal}(A):={1\over[K:{\bf Q}]}\log(\#(\omega_{\cal A}/\alpha\cdot\omega_{\cal A}))-
{1\over 2[K:{\bf Q}]}\sum_{\sigma:K\to\mC}\log|{1\over(2\pi)^{d}}
\int_{{\cal A}(\mC)_{\sigma}}
\alpha_\sigma\wedge\bar{\alpha}_\sigma|.
\label{defFH}
\end{equation}

It turns out that $h_{\rm Fal}(A)$ only depends on $A_{\bar K}$ (in particular, it does not depend on $\alpha$), so the modular height is an invariant of abelian varieties over $\bar\mQ$. The modular height plays a key role in Faltings' first proof of 
the Mordell conjecture, and it gives rise to height functions on various moduli spaces of abelian varieties. 

Now suppose in addition that we are given 
a CM field $E$, which is of degree $2d$ over $\bf Q$, and an 
embedding of rings $\phi:{\cal O}_{E}\to {\rm End}(A)$ of ${\cal O}_{E}$ into the 
endomorphism ring of $A$ (over $K$). This is equivalent to 
saying that $A$ has complex multiplications  by ${\cal O}_{E}$ 
(cf. \cite{Lang-CM}). In other words, $A$ has complex multiplication by the maximal 
order of $E$ (one can also consider complex multiplication by other orders, but we do not consider this situation here). We recall that a field is CM if it is a quadratic imaginary extension of 
a totally real number field. Such a field is endowed with an involution $c_E$, which is uniquely characterised by the fact that it extends to complex conjugation after any embedding of 
the field into $\mC.$

Suppose as well that all the embeddings $E\to \mC$ can be written as compositions 
$E\to K\to\mC$ of embeddings. We also suppose for simplicity of exposition 
that $E\subseteq K\subseteq\mC$ (ie we choose fixed embeddings of $E\to K$ and of 
$K\to\mC$) and that both $E$ and $K$ are Galois extension of $\mQ$. 

Let $G:=\Gal(E|\mQ).$ Let 
$$
\Phi:=\{\tau\in \Gal(E|\mQ)\,|\,\{t\in \Lie(A):a(t)=\tau(a)t\ \forall a\in
{\cal O}_E\}\not =0\}.
$$ 
Here $\Lie(A)$ is the Lie algebra of $A$. 
The set $\Phi$ is the {\it type} of the CM abelian 
variety $A$.  We identify 
$\Phi$ with its characteristic function $G\to\{0,1\}$ and we define $\Phi^{\vee}:G\to\{0,1\}$ by the formula $\Phi^{\vee}(\tau):=
\Phi(\tau^{-1})$. 

If $f,g:G\to \mC$ are two complex valued functions on $G$, we define
$$
\langle f,g\rangle:= {1\over \#G}\sum_{\gamma\in G}f(\gamma)\overline{g(\gamma)}
$$
($L^2$ scalar product of complex valued functions on $G$) and
$$
(f\ast g)(\lambda):={1\over \#G}\sum_{\gamma\in G} f(\gamma) g(\gamma^{-1} \lambda)
$$ 
for any $\lambda\in G$ (convolution product of complex valued functions on $G$). 

Colmez now conjectured the following:

\begin{conj}[Colmez]
The identity 
$$
{1\over d}h_{\rm Fal}(A)=-\sum_{\chi\ {\rm odd}}<\Phi*\Phi^{\vee},\chi>
[2{L'(\chi,0)\over L(\chi,0)}+
\log(f_{\chi})]
$$
holds. 
\label{conjC}
\end{conj}
Here the sum runs over all the odd irreducible Artin characters $\chi$ of $G$ (recall that 
an Artin character $\chi$ is odd iff $\chi(\gamma\circ c\circ \gamma^{-1}\circ\tau)=-\chi(\tau)$ for all
$\tau,\gamma\in G$). The 
symbol $f_{\chi}$ refers to the conductor of $\chi$ and the symbol $L(\chi,s)$ to the (Artin) $L$-function associated with $\chi$. 

Colmez proved in \cite{Colmez-Periodes} (see also \cite{Colmez-Fal}) that Conjecture \ref{conjC} holds if $G=\Gal(E|\mQ)$ is an abelian group, in the situation when $E$ is not too ramified over $2$. Additional work by Obus 
(see \cite{Obus-Col}) showed that this condition was superfluous, so the conjecture holds 
if $G$ is abelian. Note that in that case, the Artin characters in \ref{conjC} become Dirichlet characters. The conjecture of Colmez is still completely open in the situation where $E$ is not abelian. 

The refinement of Conjecture \ref{conjC} that we will soon describe is motivated 
by the fact that in the definition of the Faltings height, one considers a sum over 
all the embeddings of $K$ into $\mC$. One is thus led to ask whether there 
is a formula for each individual term in the sum (up to the natural indetermination of $\alpha$). The fact that it should be possible to make such a conjecture is 
suggested by the Arakelov theoretic proof of a weak form of the result of Colmez-Obus (ie Conjecture \ref{conjC} when $G$ is abelian) given in \cite{KR-IV}. In \cite{KR-IV}, 
the final identity is obtained by averaging a function defined over $\Gal(K|\mQ)$. 
This function does not identify each term in the sum appearing in the Faltings height, but it points to the conjectural refinement of \ref{conjC} that we shall now present. See the explanations before Remark \ref{rem1} below for more details. 

\begin{conj} There exists an integer $m\geq 1$ and an element $\beta=\beta(A,m,\phi)\in \omega_{\cal A}^{\otimes m}$
such that 

{\rm (a)} the element $\beta$ generates $\omega_{\cal A}^{\otimes m}$ (in particular, $\omega_{\cal A}^{\otimes m}\simeq\CO_K$);

{\rm (b)} for any $\sigma\in\Gal(K|\mQ)$ and $\alpha\in \omega_{\cal A}\without\{0\}$, we have 
\begin{equation}
\big|(2\pi)^{-d}\int_{{\cal A}_\sigma(\mC)} \alpha_\sigma\wedge\bar{\alpha}_\sigma\big|^{1/2}=|\sigma(\alpha^{\otimes m}/\beta)|^{1/m}\exp\Big(d\sum_{\chi\ {\rm odd}}<\Phi*\Phi^{\vee},\chi>
[2{L'(\chi,0)\over L(\chi,0)}+
\log(f_{\chi})]\Big).
\label{mainEq}
\end{equation}
\label{mainC}
\end{conj}
Here $\alpha^{\otimes m}/\beta$ is the unique element $t\in \CO_K$ such that 
$t\cdot\beta=\alpha^{\otimes m}.$ Note that if an abelian variety defined over a number field 
has complex multiplications, then it has good reduction at any place of the number field where it has semistable reduction  (see \cite{Serre-Tate-Good}). In particular, with our assumptions,  $A$ has good reduction at all the finite places of $K$.

We chose to formulate Conjecture \ref{mainC} without using the language of Arakelov theory, but it can be expressed more succinctly in that language. For the reader who knows some Arakelov theoretic terminology, the conjecture can be written as follows. The $\CO_K$-module $\omega_A$ 
can be made into a hermitian $\CO_K$-module $\bar\omega_{\cal A}$ by endowing $\omega_{\cal A}\otimes_{\CO_K,\sigma}\mC$ with the unique hermitian metric $h_\sigma(\cdot,\cdot)$ such that 
$$
h_\sigma(\eta,\eta)=|{1\over(2\pi)^{d}}
\int_{{\cal A}_\sigma(\mC)}
\eta\wedge\bar{\eta}|
$$
for any $\eta\in \omega_{\cal A}\otimes_{\CO_K,\sigma}\mC$. 
Let $\bar M$ be the hermitian $\mZ$-module, which has $\mZ$ as underlying $\mZ$-module, 
and such that $M\otimes_{\mZ}\mC$ is endowed with the hermitian metric such that 
the norm of $1\otimes 1$ is $$\exp\Big(d\sum_{\chi\ {\rm odd}}<\Phi*\Phi^{\vee},\chi>
[2{L'(\chi,0)\over L(\chi,0)}+
\log(f_{\chi})]\Big).$$

Let $\pi:\Spec(\CO_K)\to\Spec(\mZ)$ be the natural morphism. 
Then Conjecture \ref{mainC} says that there exists an integer $m\geq 1$ such that there is an isomorphism of hermitian $\CO_K$-modules $\bar\omega_{\cal A}^{\otimes m}\simeq\pi^*(\bar M^{\otimes m}).$ 

In \cite[Example 7.1]{MR-Biel25}, a weak form of this last statement is proven in the situation where $G$ is abelian. The proof is based on the Lefschetz fixed formula in Arakelov theory 
proven in \cite{KR-I}. More precisely, in the terminology of \cite[Def. 5.3 and after]{MR-Biel25}, it is shown in \cite[Example 7.1]{MR-Biel25}  that if $G$ is abelian, then the first arithmetic Chern class of $\bar\omega_{\cal A}$ in the 
group $\widehat{\rm CH}^1_{\bar\mQ}(\CO_{K,2{\rm disc}(E)})$ is equal 
to the class of the constant $$-d\sum_{\chi\ {\rm odd}}2<\Phi*\Phi^{\vee},\chi>
{L'(\chi,0)\over L(\chi,0)}.$$ Here ${\rm disc}(E)$ is the discriminant of $E$. 
This implies in particular that the first arithmetic Chern classes of $\bar\omega_{\cal A}$ and $\pi^*(\bar M)$ in $\widehat{\rm CH}^1_{\bar\mQ}(\CO_{K,2{\rm disc}(E)})$ coincide. However, this statement is not strong 
enough to imply Conjecture \ref{mainC} when  $G$ is abelian. 

\begin{rem}\rm 
\label{rem1}
(1) Both sides of the equality in (b) are multiplied by $|\sigma(a)|$ if we multiply 
$\alpha$ by $a\in\CO_K$. So if (b) holds for a given $m$ and $\beta$ and some choice of $\alpha$, then it holds for $m$, $\beta$ and any other choice of $\alpha.$

(2) There is always an integer $m\geq 1$, such that 
$\omega_{\cal A}^{\otimes m}\simeq\CO_K$, because the class group of 
$\CO_K$ is finite. However, the integer $m$ required by Conjecture \ref{mainC} 
might be larger than the exponent of the class group.

(3) If $m$ is given, then the generator $\beta$ (if it exists) is uniquely determined 
up to multiplication by a root of unity lying in $K$. Indeed, suppose that 
$\beta$ and $\beta'$ are two generators of $\omega_{\cal A}^{\otimes m}$ satisfying \refeq{mainEq}. 
Let $\beta'=u\cdot \beta$, so that $u$ is a unit of $\CO_K$. 
From \refeq{mainEq}, we see that $|\sigma(u)|=1$ for all 
$\sigma:K\to\mC$. Hence $u$ is a root of unity by the Dirichlet unit theorem.

(4) If (a) and (b) hold for some $m$ and $\beta$ and for some 
$m'$ and $\beta'$, then the conjecture holds for 
$m+m'$ and $\beta\otimes\beta'$. In particular, if (a) and (b) hold for some $m$ and 
$\beta$, then (a) and (b) hold for $mt$ and $\beta^{\otimes t}$. Hence, if there exists an integer $m$, which is minimal 
with the property that there exists a $\beta$ such that (a) and (b) hold, then 
the other integers $m'$ with that property are all multiples of $m$.

(5) The Arakelov theoretic approach to the formula of Colmez suggests that there should be an explicit formula for an integer $m$ such that Conjecture \ref{mainC} holds for $A$ and $\phi$. We are not able to make an educated guess for such a formula yet. Here is a very wild guess (which we do not dare to state as a conjecture) about the shape of such an $m$. 
Let $h^{-}_E$ be the relative class number of $E$. 
Suppose that there is an integer $m_0=m_0(A,\phi)$, which is minimal with the property that Conjecture \ref{mainC} holds for $A,\phi$ and $m=m_0.$ Then one might speculate that if $p|m_0$ and $p>2$, then either $p$ divides the discriminant of $E$ or $p$ divides $[E:\mQ]\cdot h^{-}_E.$ 
We will not offer any motivation for this conjecture, which is based on very heuristic considerations. Proposition \ref{propELL} below shows that it holds if $p>3$ and $A$ is an elliptic curve.   

\end{rem} 

\begin{lemma} Conjecture \ref{mainC} implies Conjecture \ref{conjC}.\end{lemma}
\beginProof From the definitions, we see that it is sufficient to show that 
\begin{equation}
\prod_{\sigma:K\to\mC}|\sigma(\alpha^{\otimes m}/\beta)|^{1/m}=\#(\omega_{\cal A}/\alpha.\omega_{\cal A}).
\label{midEq2}
\end{equation}
To show that \refeq{midEq2} holds, note first that 
\begin{equation}
\#(\omega_\CA^{\otimes m}/\alpha^{\otimes m}\cdot\omega_\CA^{\otimes m})=\#(\omega_{\cal A}/\alpha\cdot\omega_{\cal A})^m.
\label{midEq}
\end{equation}
To establish \refeq{midEq}, note that $\omega_{\cal A}/\alpha.\omega_{\cal A}$ 
is isomorphic as a $\CO_K$-module to a module of the form 
$\oplus_{\Fp\in{\rm Spec}(\CO_K)}\CO_K/\Fp^{n_\Fp}$, so that 
$$
\omega_\CA^{\otimes m}/\alpha^{\otimes m}\cdot\omega_{\cal A}^{\otimes m}
\simeq \oplus_{\Fp\in{\rm Spec}(\CO_K)}\CO_K/\Fp^{m n_\Fp}.
$$
The fact that $\#(\CO_K/\Fp^n)=\#(\CO_K/\Fp)^n$ holds for any $n\geq 1$ now immediately implies  that \refeq{midEq} holds. To prove \refeq{midEq2}, we are thus reduced to showing that 
$$
\prod_{\sigma:K\to\mC}|\sigma(\alpha^{\otimes m}/\beta)|=\#(\omega_\CA^{\otimes m}/\alpha^{\otimes m}.\omega_\CA^{\otimes m})
$$
and since $\omega_{\cal A}^{\otimes m}\simeq\CO_K$, this amounts to showing 
that 
$$
\prod_{\sigma:K\to\mC}|\sigma(a)|=\#(\CO_K/a\CO_K)
$$
for any $a\in\CO_K$. But this follows from standard properties of the ideal norm.\endProof

\section{Siegel modular forms}
\label{secSIEG}

One interesting aspect of Conjecture \ref{mainC} is that it implies a reciprocity law for Siegel modular forms evaluated at CM points in Siegel upper half space. Conversely, if one can exhibit a Siegel modular which satisfies a reciprocity law at a CM point and which has suitable integrality properties, then one can use it to show that Colmez's conjecture implies Conjecture \ref{mainC} for that CM point. We shall illustrate this principle for elliptic curves in section \ref{secELL}, where it will be shown that the discriminant modular form has all the required properties. 

We first recall some classical facts. Let $f$ be a holomorphic Siegel modular  form  of 
weight $w$ for ${\rm Sp}(2g,\mZ)$. The integer $g$ is called the {\it genus} of $g$. Suppose that the coefficients of the $q$-expansion $f$ all lie in 
$R$, where $R\subseteq\mC$ is a subring. If 
$\CB$ is a principally polarised abelian abelian scheme over $R$, the modular form $f$ then gives rise to a section $f(\CB)$ of the sheaf $(\Omega^g_\CB)^{\otimes w}$ over $\CB$.  
This follows from (deep) results of Chai and Faltings - see \cite[in particular p. 141]{CF-Deg}. 
One should not confuse $f(\CB)$ 
with the evaluation $f(\Tau)$ of $f$ at an element $\Tau$ of the Siegel upper half space $\CH_g$, which is a complex number (but this should not lead to any confusion). If $R=\mC$, we shall write $\langle \cdot,\cdot\rangle_\Pet=\langle \cdot,\cdot\rangle_{w,\CB,\Pet}$ for the metric on the $\mC$-vector space of global sections of $(\Omega^g_{\CB})^{\otimes w}$ which is such that 
 $$
\langle \eta^{\otimes w},\eta^{\otimes w}\rangle_\Pet=|(2\pi)^{-g}\int_{\CB(\mC)} \eta\wedge\bar{\eta}\big|^w
$$
for any global section $\eta$ of $\Omega^g_{\CB}$. We shall then write 
$\lVert\cdot\rVert_\Pet$ for the associated norm. If $\Tau\in\CH_g$ is a point of Siegel half space representing $\CB$, the formula  
$$
\lVert f(\CB)\rVert_\Pet=|f(\Tau)|\det(\Im(\Tau))^{w/2}(4\pi)^{gw/2}.
$$
then holds. See \cite[before Lemma 2.1]{Yang-CS} for lack of a more standard reference. 

\begin{prop}
Let the terminology of Conjecture \ref{mainC} hold.  Let $M\subseteq K$ be a number field and let $S\subseteq\CO_M\backslash\{0\}$ be a multiplicative set. Let $f$ be a Siegel modular 
form of weight $m$ and genus $d$, with coefficients  in $\CO_{M,S}$. For all $\sigma\in\Gal(K|M)$ let $\Tau_\sigma\in\CH_d$ be a representative of $A_\sigma(\mC)$ in Siegel upper half space.

Suppose that 
$\CA$ is endowed with a principal polarisation and suppose that $f(\CA_{\CO_{K,S}})$ generates $\omega_{\CA_{\CO_{K,S}}}^{\otimes m}$. 

Then, if Conjecture \ref{mainC} holds, there exists an element $\rho=\rho(f,A,E)\in U(\CO_{K,S})$ such that
\begin{eqnarray}
|\sigma(\rho)|&=&{|f(\Tau_\sigma)|\det(\Im(\Tau_\sigma))^{m/2}(4\pi)^{dm/2}\over \exp\Big(m d\sum_{\chi\ {\rm odd}}<\Phi*\Phi^{\vee},\chi>
[2{L'(\chi,0)\over L(\chi,0)}+
\log(f_{\chi})]\Big)}\nonumber\nonumber\\&&\nonumber\\&=:&{\lVert f(A_\sigma)\rVert_\Pet\over \exp\Big(m d\sum_{\chi\ {\rm odd}}<\Phi*\Phi^{\vee},\chi>
[2{L'(\chi,0)\over L(\chi,0)}+
\log(f_{\chi})]\Big)}
\label{eqREC}
\end{eqnarray}
for all $\sigma\in\Gal(K|M).$ 
\label{propRL}
\end{prop}

Here $U(\bullet)$ is the group of units of the commutative ring $\bullet$. Note also that 
any principal polarisation of $A$ over $M$ extends uniquely to a principal polarisation of $\CA$ over $\CO_{K,S}$. 

\begin{rem}\rm (1) Note that if $S:=\{1\}$ and $M=\mQ$, then the element $\rho$ is determined up to multiplication by a root of unity in $K$ by $f$, $\phi$ and  by the condition \refeq{eqREC}. 

(2) The condition that $f(A)$ generates $\omega_{\CA_{\CO_{K,S}}}^{\otimes m}$ can be understood 
as an integrality condition. Denote by $\CA_{d}$ the algebraic 
stack classifying  principally polarised abelian schemes of relative dimension $d$. Then the condition that  $f(A)$ generates $\omega_{\CA_{\CO_{K,S}}}^{\otimes m}$ amounts to asking whether the element of 
$\CA_g(K)$ defined by $A$ is $S$-integral with respect to the divisor defined by $f$. 
When $A$ is an elliptic curve and $f$ is the (rescaled) modular discriminant, this 
always holds because the divisor $f$ is empty. This will be exploited in section \ref{secELL}. It is not clear to us whether in the situation of Conjecture 
\ref{mainC}, this integrality condition can always be achieved. Note that when 
$d>1$, then the divisor defined by $f$ is relatively ample on a suitable \'etale covering of 
$\CA_g$ (this is the main result of 
\cite{MB-Pinceaux}), so that integrality becomes a thorny issue if $d>1$.
\end{rem}

\beginProof (of Proposition \ref{propRL}). The proof is straightforward. Conjecture \ref{mainC} implies that 
$$
\lVert\beta_\sigma\rVert=\exp\Big(m d\sum_{\chi\ {\rm odd}}<\Phi*\Phi^{\vee},\chi>
[2{L'(\chi,0)\over L(\chi,0)}+
\log(f_{\chi})]\Big)
$$
for all $\sigma\in\Gal(K|M).$ Now define $\rho\in\CO_{K,S}$ by the formula 
$$
\rho\cdot \beta=f(\CA).
$$
We then have
$$
\lVert\sigma(\rho)\cdot\beta_\sigma\rVert_\Pet=|\sigma(\rho)|\cdot \lVert\beta_\sigma\rVert_\Pet=\lVert f(\CA)_\sigma\rVert_\Pet=
\lVert f(\CA_\sigma)\rVert_\Pet
$$
from which \refeq{eqREC} follows.
\endProof

The following reciprocity law at CM points now readily follows for the Siegel modular form $f$:

\begin{cor}
Let the assumptions and terminology of Proposition \ref{propRL} hold and suppose in addition that $S$ is invariant under conjugation. Then if Conjecture \ref{mainC} holds, 
there are elements $\{a_\sigma\in \mC\}_{\sigma\in\Gal(K|M)}$ with the following properties:

{\rm (a)} $|a_\sigma|=1$ for all $\sigma\in\Gal(K|M)$;

{\rm (b)} $a_\sigma^2\in U(\CO_{K,S})$ for all $\sigma\in\Gal(K|M)$;

{\rm (c)} the elements
$$
u(\sigma):= a_\sigma{|f(\Tau_\sigma)|\det(\Im(\Tau_\sigma))^{m/2}\over 
|f(\Tau_\Id)|\det(\Im(\Tau_\Id))^{m/2}}
$$
lie in $U(\CO_{K,S})$ for all $\sigma\in\Gal(K|M)$ and satisfy the reciprocity law
$$
\tau(u(\sigma))=u({\tau\sigma})/u(\tau)
$$
for all $\sigma,\tau \in\Gal(K|M).$
\label{corSRL}
\end{cor}

\beginProof In the notation of Proposition \ref{propRL}, let $u(\sigma):=\sigma(\rho)/\rho$ and $a_\sigma:=u(\sigma)/|u(\sigma)|$.\endProof

When $d=1$, $M=E$ and $S=\CO_K\backslash\{0\}$, Corollary \ref{corSRL} follows from Shimura's reciprocity law (see \cite[chap. 11, par. 2, Th. 6]{Lang-Ell}). When $S=\{1\}$, $M=E$, $d=1$ and $f$ is the discriminant modular form, Corollary \ref{corSRL} follows from work of Hasse - see Proposition \ref{propHDS} below. 

\begin{rem}\rm (1) It should be possible to approach the conclusion of Corollary \ref{corSRL} using generalisations of Shimura's reciprocity law. If $M=\mQ$, such reciprocity laws would have to be based on the general version of the fundamental theorem of complex multiplication (for all automorphisms of $\mC$, see \cite{Schappacher-CM}), and not on the older version (due to Shimura and Taniyama) 
which assumes that the automorphism of $\mC$ fixes the reflex field of $E$ with respect to 
the type $\Phi$. In \cite[par 3.2, Lemma 4]{GdS-Theta}, Shimura's reciprocity law is used to prove part of the conclusion of Corollary \ref{corSRL} in the situation where $E$ is a quartic CM field and $f$ is Igusa's modular form $\chi_{10}(\cdot)$. 
We hope to come back to this line of investigation in a subsequent article.

(2) Shimura's monomial period relations imply that the quantities ${|f(\Tau_\sigma)|\det(\Im(\Tau_\sigma))^{m/2}\over 
|f(\Tau_\Id)|\det(\Im(\Tau_\Id))^{m/2}}$ are all algebraic (we skip the details), as implied by 
Corollary \ref{corSRL} (under the assumption that Conjecture \ref{mainC} holds). 
We are grateful to Pierre Colmez for pointing this out to us. 
Pierre Colmez also suggested that if one used the refined monomial relations established in his recent article \cite[Th. 12]{Colmez-AbsUn} (which 
are derived from the existence of $p$-adic analogs of Deligne's results \cite{DMOS} on absolute Hodge cycles), it might be possible to show that 
the same quantities actually lie in $|U(\CO_{K,S})|$, as predicted by Corollary \ref{corSRL}. Again, we hope to get back to this in a later article.
 \end{rem}

The following Lemma is a partial converse of Corollary \ref{corSRL}. 
We will apply this lemma to prove Conjecture \ref{mainC} (and in particular Corollary \ref{corSRL}) for CM elliptic curves in 
the next section. 

\begin{lemma} Let the terminology of Conjecture \ref{mainC} hold.  Let $M\subseteq K$ be a number field and let $S\subseteq\CO_M\backslash\{0\}$ be a multiplicative set. Let $f$ be a Siegel modular 
form of weight $w$ and genus $d$, with coefficients  in $\CO_{M,S}$. Suppose that 
$\CA$ is endowed with a principal polarisation and suppose that $f(\CA_{\CO_{K,S}})$ generates $\omega_{\CA_{\CO_{K,S}}}^{\otimes w}$. For all $\sigma\in\Gal(K|M)$ let $\Tau_\sigma\in\CH_d$ be a representative of $A_\sigma(\mC)$ in Siegel upper half space.  

Finally, suppose that there exist elements $\{u(\sigma)\in U(\CO_{K,S})\}_{\sigma\in\Gal(K|M)}$ such that
\begin{equation}
|u(\sigma)|={|f(\Tau_\sigma)|\det(\Im(\Tau_\sigma))^{w/2}\over 
|f(\Tau_\Id)|\det(\Im(\Tau_\Id))^{w/2}}\label{eqREC2}.
\end{equation}
and such that 
$$
\tau(u(\sigma))=u({\tau\sigma})/u(\tau)
$$
for all $\sigma,\tau\in\Gal(K|M)$. 

Let $e=e(u(\bullet))$ be a multiple of the order in $H^1(\Gal(K|M),U(\CO_{K,S}))$ of the cocycle 
with values in $U(\CO_{K,S})$ defined by $u{(\bullet)}$. 

Then:

{\rm (a)} There exists a generating element $\beta\in\omega_{\CA_{\CO_{K,S}}}^{\otimes we}$ such that 
$\lVert\beta\rVert_\Pet=\lVert\beta_\tau\rVert_\Pet$ for all $\tau\in\Gal(K|M).$

{\rm (b)} Suppose that $S=\{1\}.$ If $M|\mQ$ is an imaginary quadratic field or if $M=\mQ$, then Conjecture \ref{mainC} holds for $\beta$ and $m=we$.
\label{lemAUTR}
\end{lemma}
\beginProof 
By construction, there is an element $u\in U(\CO_{K,S})$ such that 
$u(\sigma)^e=\sigma(u)/u$ for all $\sigma\in\Gal(K|M).$ Such an element 
can actually easily be chosen explicitly. Eg, if $e=\#\Gal(K|M)$, the element 
$$
u:=\big(\prod_{\tau\in \Gal(K|M)} u(\tau))^{-1}
$$
has the required properties, as shown by the computation
$$
\sigma(u)=\big(\prod_{\tau\in\Gal(K|M)} \sigma(u(\tau))\big)^{-1}=
\big(\prod_{\tau\in\Gal(K|M)} {u(\sigma\tau)\over u(\sigma)}\big)^{-1}=
\big(\prod_{\tau\in\Gal(K|M)} {u(\tau)\over u(\sigma)}\big)^{-1}=u\cdot u(\sigma)^{\#\Gal(K|M)}.
$$
We let
$$
\beta:=f(\CA_{\CO_{K,S}})^{\otimes e}\cdot u^{-1}=f^e(\CA_{\CO_{K,S}})\cdot u^{-1}\in\omega_{\CA_{\CO_{K,S}}}^{\otimes we}.
$$
The element $\beta$ generates  $\omega^{\otimes we}_{\CA_{\CO_{K,S}}}$, because $u$ is a unit of $\CO_{K,S}$ and because by assumption $f^{e}$ generates $\omega^{\otimes we}_{\CA_{\CO_{K,S}}}$. 

Let $\tau\in\Gal(K|M).$ We compute
\begin{eqnarray*}
\lVert\beta_\tau\rVert_\Pet&=&|f^e(\Tau_\tau)|\det(\Im(\Tau_\tau))^{ew/2}(4\pi)^{dew/2}|\tau(u^{-1})|\\&=&
|f^e(\Tau_\tau)|\det(\Im(\Tau_\tau))^{ew/2}(4\pi)^{dew/2}|u^{-1}||u/\tau(u)|\\
&=&|f^e(\Tau_\Id)|\det(\Im(\Tau_\Id))^{ew/2}(4\pi)^{dew/2}|u^{-1}|=:\lVert\beta\rVert_\Pet.
\end{eqnarray*}
This establishes (a). 

Now suppose that $M=\mQ$ and $S=\{1\}.$ Then by (a) we have
\begin{eqnarray*}
\lVert\beta\rVert_\Pet&=&\big(\prod_{\sigma\in\Gal(K|\mQ)}\lVert \beta_\sigma\rVert_\Pet)\big)^{1/[K:\mQ]}=\exp(h_{\rm Fal}(E))^{we}
\\&=&\exp\Big(we d\sum_{\chi\ {\rm odd}}<\Phi*\Phi^{\vee},\chi>
[2{L'(\chi,0)\over L(\chi,0)}+
\log(f_{\chi})]\Big), 
\end{eqnarray*}
so $\beta$ has all the properties required by Conjecture \ref{mainC} if $m=we.$ 

Now suppose that $M$ is an imaginary quadratic field and $S=\{1\}.$ Note first that the Peterson metric is invariant under complex conjugation (this can be seen from the formula). 
Denoting by $c$ the restriction of complex conjugation to $K$, we may thus compute
\begin{eqnarray*}
\lVert\beta\rVert_\Pet&=&\big(\prod_{\sigma\in\Gal(K|M)}\lVert \beta_\sigma\rVert_\Pet)\big)^{1/(2[K:M])}\big(\prod_{\sigma\in\Gal(K|M)}\lVert \beta_{c\sigma}\rVert_\Pet)\big)^{1/(2[K:M])}\\&=& 
\big(\prod_{\sigma\in\Gal(K|\mQ)}\lVert \beta_\sigma\rVert_\Pet)\big)^{1/[K:\mQ]}=\exp(h_{\rm Fal}(E))^{we}
\\&=&\exp\Big(we d\sum_{\chi\ {\rm odd}}<\Phi*\Phi^{\vee},\chi>
[2{L'(\chi,0)\over L(\chi,0)}+
\log(f_{\chi})]\Big), 
\end{eqnarray*}
and if $\tau\in\Gal(K|M)$, we have $\lVert\beta\rVert_\Pet=\lVert\beta_\tau\rVert_\Pet=
\lVert\beta_{c\tau}\rVert_\Pet$ from (a) and the above remark, so that $\lVert\beta\rVert_\Pet=\lVert\beta_\sigma\rVert_\Pet$ 
for all $\sigma\in\Gal(M|\mQ).$ 
So again $\beta$ has all the properties required by Conjecture \ref{mainC} if $m=we.$ 
\endProof

\begin{rem}\rm (1) The  construction of the unit $u$ in the proof of Lemma \ref{lemAUTR} amounts to solving a "descent" problem for 
the hermitian $\CO_K$-module $\bar\omega^{\otimes we}_{\CO_{K,S}}$ from $\CO_K$ to 
$\mZ$. The element $f(\CA_{\CO_{K,S}})^{\otimes e}$ does not necessarily 
descend to $\mZ$, but the element $f(\CA_{\CO_{K,S}})^{\otimes e}\cdot u^{-1}$ does.

(2) In all of the above, one can consider modular forms 
for arithmetic groups other than ${\rm Sp}(2g,\mZ)$ if one can construct good integral 
models (possibly as algebraic stacks) of the corresponding Shimura varieties, and also if there is a way to determine which modular forms come from the model (in this case it can be checked by looking at the ring of definition of the coefficients of the $q$-expansion). For example, Moret-Bailly has shown in \cite[Th. 0.5]{MB-Sur} 
that the eighth power of the Riemann theta function, which is modular with respect to the group $\Gamma(1,2)\subseteq {\rm Sp}(2g,\mZ)$ of Igusa, has the required properties with respect to a natural integral 
model of the corresponding Shimura variety, so Proposition \ref{propRL}, Corollary \ref{corSRL} and Lemma \ref{lemAUTR} can be formulated 
for the Riemann theta function.\end{rem}

\section{The case of elliptic curves}
\label{secELL}

We keep the terminology of Conjecture \ref{mainC}.

\begin{prop} Conjecture \ref{mainC} holds if $d=1$ and $m=12h^2$, where 
$h$ is the class number of $E$. 
\label{propELL} \end{prop}

Before we move to the proof of Proposition \ref{propELL}, we shall introduce some terminology and recall some classical results of the theory of 
complex multiplication of elliptic curves. Since $d=1$, we have $E=\mQ(\sqrt{-d_0})$, where 
$d_0$ is a square free positive integer. Let $-D$ be the discriminant of the number field $E$. It is classical 
that $D=d_0$ if $d_0=-1\,(\md{4})$, and $D=4d_0$ otherwise. 
We also 
recall that one can show that $\CO_E=\mZ[(1+i\sqrt{d_0})/2]$ if $d_0=-1\,(\md{4})$ and 
$\CO_E=\mZ[i\sqrt{d_0}]$ otherwise. 
Let $\Cl(E)$ be the class group of $\CO_E$ and let $h:=\#\Cl(E)$ be the class number of $E$. Also, let $w_E$ be the cardinality of 
the group of roots of unity in $E$. We note that $w_E|12$. Finally, we note that 
by Dirichlet's unit theorem, all the units of $\CO_E$ are roots of unity.

Let $(2\pi)^{-12}\Delta(\cdot)=\eta^{24}(\cdot)$ be the rescaled discriminant modular form, which is a Siegel modular 
form of weight $12$ on the Siegel upper half space $\CH_1$. The Fourier coefficients (coefficients of the $q$-expansion) of 
$(2\pi)^{-12}\Delta(\cdot)$  are all integers. Also, if $B$ is an elliptic curve with good reduction at all places over a number field $L$, the element 
$\Delta(B)\in\omega^{\otimes 12}_B$ arises from a generating element of 
$\omega^{\otimes 12}_{\CB}$, where $\CB$ is the N\'eron model of $B$ over 
$\CB$.  

If $\Lambda\subseteq\mC$ is a lattice with positively oriented basis $z_1,z_2$, and 
$f$ is a modular form of genus $1$ and of weight $w$, we shall 
write 
\begin{equation}
f(\Lambda):=z_2^{-w}f(z_1/z_2)
\label{eqDELH}
\end{equation}
where by construction $z_1/z_2$ is in the upper half plane. This ("homogenous") definition of $\Delta(\Lambda)$ does not depend on the choice of 
the positively oriented basis $z_1,z_2$ of $\Lambda$. 

Let now $\mFa$ be a fractional ideal in $E$. Let $b(\mFa)$ be a generator 
of $\mFa^{h}$ as a $\CO_E$-module. We define
$$
\rho(\mFa):={\Delta(\CO_E)^{h}\over b(\mFa)^{12}\,\Delta(\mFa)^{h}}.
$$
Note that $\beta(\mFa)^{12}$ does not depend on the choice of $\beta(\mFa)$ 
because $w_E|12$ (note that two different choices of $\beta(\mFa)$ differ by a unit, and hence a root of unity in our situation). 
Furthermore, note that $\rho(\mFa)$ only depends on the ideal class of $\mFa$, because 
$\beta((a)\mFa)^{12}\,\Delta((a)\mFa)^{h}=\beta(\mFa)^{12}\,\Delta(\mFa)^{h}$ by \refeq{eqDELH}.

We also recall that $E(j(A))=:H$ is the Hilbert class field of $E$. Here $j(\cdot)$ is the $j$-invariant, a modular form of weight $0$. A consequence of this is that the field of definition $K$ of $A$ contains $H.$ The Artin reciprocity map $(\bullet,H/E):\Cl(E)\to\Gal(H|E)$ gives an isomorphism 
between the class group of $E$ and the Galois group of $H$ over $E$. 
We have $A_\Id(\mC)\simeq\mC/\mFb_A$ for some ideal $\mFb_A$ of $\CO_E$, and the formula $A_\sigma(\mC)\simeq\mC/(\sigma|_H,H/E)^{-1}\mFb_A$ holds for all 
$\sigma\in\Gal(K|E).$ In particular, if $f$ is a modular form of genus $1$ and of weight $0$, we have $f(\Tau_\sigma)=f((\sigma|_H,H/E)^{-1}\mFb_A)$ (where as usual $\Tau_\sigma\in\CH_1$ represents $A_\sigma$). Hence we have $j(A)=j(\mFb_A)$. 

Finally, note that since to prove Proposition \ref{propELL}, we can replace $A$ by one of its twists without restriction of generality, we may therefore assume that $\mFb_A=\CO_E$. 

{\it So we shall assume until the end of the section that $A_\Id(\mC)\simeq\mC/\CO_E$.}

See \cite[chap. 10]{Lang-Ell} for more details on all the material described above.

We can now state the

\begin{prop}[Hasse, Deuring, Siegel]

Let $\mFa$ be a fractional ideal of $\CO_E$. 

{\rm (1)} We have $\rho(\mFa)\in U(\CO_H).$

{\rm (2)} For any  fractional ideal $\mFb$ of $E$, we have 
$$
\rho(\mFa)^{(\mFb,K)}=
\rho(\mFb^{-1}\mFa)\rho(\mFb^{-1})^{-1}.
$$

{\rm (3)} The formula
$$
|\rho(\mFa^{-1})|={\det(\Im(\Tau_\Id))^{6h}|\Delta(\Tau_\Id)^h|\over \det(\Im(\Tau_{(\mFa,H/E)}))^{6h}|\Delta(\Tau_{(\mFa,H/E)}))^h|}
$$
holds.
\label{propHDS}
\end{prop}
\beginProof For (2) see \cite[Th. 12.1.1, p. 161]{Lang-Ell} and for (1), see  \cite[Chap. II, sec. 2, in particular p. 79]{Siegel-Adv}. We prove (3). Note that 
by the discussion at the beginning of this section, $\CO_E$ is generated as a lattice by $z_E$ and $1$, where $\Im(z_E)>0.$ Also, since $\rho(\mFa^{-1})$ only depends on the ideal class of $\mFa^{-1}$, we may divide  $\mFa^{-1}$ by an appropriate element of $E$ and also suppose that $\mFa^{-1}$ is generated as a lattice by 1 and $z_{\mFa^{-1}}$, say, where $\Im(z_{\mFa^{-1}})>0.$  In particular, we then have $\CO_E\subseteq\mFa^{-1}.$ We have to show that 
$$
|b(\mFa^{-1})|^2=(\Im(z_{\mFa^{-1}})/\Im(z_E))^h.
$$
Since the Galois group of $E$ over $\mC$ is given by complex conjugation, we 
have $|b(\mFa^{-1})|^2=N_{E|\mQ}(\mFa^{-h})=N_{E|\mQ}(\mFa^{-1})^h$, where
$N_{E|\mQ}(\cdot)$ is the absolute ideal norm. Hence we have 
to show that 
$$
N_{E|\mQ}(\mFa^{-1}):=\#(\mFa^{-1}/\CO_E)^{-1}=\Im(z_{\mFa^{-1}})/\Im(z_E)
$$
Now, if $M:\mFa^{-1}\to\mFa^{-1}$ is map of $\mZ$-modules such that 
$M(\mFa^{-1})=\CO_E$, we have $|\det(M)|=\#(\mFa^{-1}/\CO_E)$, and thus we are reduced to show that 
$$
|\det(\begin{bmatrix} \Re(z_E) & 1 \\ \Im(z_E) & 0 \\ \end{bmatrix}\begin{bmatrix}\Re(z_{\mFa^{-1}}) & 1 \\ \Im(z_{\mFa^{-1}}) & 0 \\ \end{bmatrix}^{-1})|^{-1}=\Im(z_{\mFa^{-1}})/\Im(z_E)
$$
which is a straightforward calculation.\endProof

\beginProof (of Proposition \ref{propELL}). By Lemma \ref{lemAUTR}, 
to prove Proposition \ref{propELL}, it is sufficient to exhibit elements $u(\sigma)\in U(\CO_{K})$ such that
\begin{equation}
|u(\sigma)|={|\Delta(\Tau_\sigma)\det(\Im(\Tau_\sigma))^{6}|\over |\Delta(\Tau_\Id)\det(\Im(\Tau_\Id))^{6}|}=:{\lVert \Delta(A_\sigma)\rVert_\Pet\over \lVert \Delta(A_\Id)\rVert_\Pet}
\label{eqREC3}
\end{equation}
and such that 
$$
\tau(u(\sigma))=u({\tau\sigma})/u(\tau)
$$
for all $\sigma,\tau\in\Gal(K|E).$ 
Here, as before, $\Tau_\sigma\in\CH_d$ is a representative 
of $A_\sigma(\mC)$ in $\CH_1$ (which is just the upper half space of $\mC$ in this case).  

For each $\sigma\in\Gal(K|E)$ define
$$
u(\sigma):=\rho(\mFa(\sigma)^{-1})^{-1}
$$
where  $(\mFa(\sigma),H/E)=\sigma|_{H}$. Proposition \ref{propHDS} shows that 
all the assumptions of Lemma \ref{lemAUTR} are satisfied for $f=\eta^{24h}(\cdot)$, $w=12h$ 
and $M=E$. By construction, the 
image of $u{(\bullet)}$ in $H^1(\Gal(K|E),U(\CO_K))$ lies in the image of the composition
$$
H^1(\Gal(H|E),U(\CO_H))\to H^1(\Gal(K|E),U(\CO_H))\to H^1(\Gal(K|E),U(\CO_K))$$
where the first map is the pull-back map associated with the homomorphism 
$\Gal(K|E)\to \Gal(H|E)$ and the second map arises from the fact that goup cohomology is a functor. Finally, the order of any element of $H^1(\Gal(H|E),U(\CO_H))$ divides $h=\#\Gal(H|E)$, so we may choose $e=h$. 
Hence Conjecture \ref{mainC} holds for $12h^2$. 
\endProof

Proposition \ref{propELL} can be rewritten in terms of the $\Gamma$ function. 
Unravelling the definitions and applying a formula of Lerch (see \cite[Eq. (9), p. 4]{Duke-Kronecker}) together with a standard class number formula (see \cite[Th. 4.17 and before]{Washington-Cyc}), we obtain from Conjecture \ref{conjC} (which 
is true in this situation, since $E$ is abelian over $\mQ$) that

\begin{eqnarray}
h_{\rm Fal}(E)&=&-\langle {1\over 4}+{1\over 4}\chi,\chi\rangle\big[2{L'(\chi,0)\over L(\chi,0)}+\log(D)\big]=-{1\over 2}{L'(\chi,0)\over L(\chi,0)}-{1\over 4}\log(D)\nonumber\\&=&
-{1\over 2}\big[-\log(D)+{w_E\over 2h}\sum_{j=1}^{D-1}\legendre{-D}{j}\log(\Gamma({j\over D})\big]-{1\over 4}\log(D)\nonumber\\&=&
{1\over 4}\log(D)-{w_E\over 4h}\sum_{j=1}^{D-1}\legendre{-D}{j}\log(\Gamma({j\over D})
\label{eqGAMMA}
\end{eqnarray}

where $\legendre{\cdot}{\cdot}$ is the Kronecker symbol and $\Gamma(\cdot)$ is the classical $\Gamma$ function.

Applying the exponential function to \refeq{eqGAMMA}, one recovers the Chowla-Selberg formula (as explained in the introduction of \cite{Colmez-Periodes}):

\begin{theor}[Lerch, Chowla, Selberg \cite{CS-Epstein}] 
We have
\begin{eqnarray*}
\big[\prod_{\sigma\in\Gal(K|\mQ)}\lVert \eta^{24}(A)_\sigma\rVert_\Pet\big]^{1\over 12\cdot 2h}&=&\big[\prod_{\mFa\in\Cl(E)}\Im(z_{\mFa})|\eta(z_{\mFa})|^{4}(4\pi)\big]^{1\over 2h}=
D^{-1/4}\big[\prod_{j=1}^{D-1}\Gamma({j\over D})^{\legendre{-D}{j}}\big]^{{w_E\over 4h}}.
\end{eqnarray*}
Here $z_{\mFa}$ is an element of the upper half plane such that the lattice generated by $z_{\mFa}$ and $1$ is congruent to $\mFa$.\label{CSTHEOR}
\end{theor}
 We conclude from Proposition \ref{propELL} that there exists a unit $u\in U(\CO_K)$ such that 
 $$
 \lVert (u^{-1}\cdot \eta^{24h^2}(A))_\sigma\rVert_\Pet=D^{-3h^2}\big[\prod_{j=1}^{D-1}\Gamma({j\over D})^{\legendre{-D}{j}}\big]^{{3w_Eh}}
$$
for all $\sigma\in\Gal(K|\mQ).$ 
By the calculation presented as the beginning of Lemma \ref{lemAUTR}, the element $$
u_c:=\prod_{\mFa\in\Cl(E)}\rho(\mFa)
$$
is a possible choice of $u$, and any other choice of $u$ will differ from $u_c$ 
by a root of unity in $K$. Note that the unit $u$ was computed explicitly (up to multiplication by a root of unity in $K$) in various 
special cases in the articles \cite{CH-Eval}, \cite{VDP-W-Values}, and in \cite{HKW-CS} (see also the bibliography of these articles for further references). However, the fact that 
$$
{ \lVert \eta^{24h^2}(A)\rVert_\Pet\over D^{-3h^2}\big[\prod_{j=1}^{D-1}\Gamma({j\over D})^{\legendre{-D}{j}}\big]^{{3w_Eh}}}
$$
is a unit in general (leaving aside the reciprocity law we established) does not seem to have been mentioned before in the literature (but it was likely 
known to the experts).

Finally, note that we have $u_c\in\mR$. 
This follows from that the fact that \mbox{$\overline{\rho(\mFa)}=\rho(\bar\mFa)=
\rho(\mFa^{-1})$} for all fractional ideals $\mFa$ of $E$. We leave the proof of these identities as an exercise for the reader. So we have 
$$
\pm u_c={ \lVert \eta^{24h^2}(A)\rVert_\Pet\over D^{-3h^2}\big[\prod_{j=1}^{D-1}\Gamma({j\over D})^{\legendre{-D}{j}}\big]^{{3w_Eh}}}:={ \Im(z_E)^{6h^2}|\eta^{24}(z_E)|^{h^2}(4\pi)^{6h^2}\over D^{-3h^2}\big[\prod_{j=1}^{D-1}\Gamma({j\over D})^{\legendre{-D}{j}}\big]^{{3w_Eh}}}.
$$
Here, as before, $D=d_0$ and $z_E=(1+i\sqrt{d_0})/2$ if $d_0=-1\,(\md{4})$, and 
$D=4d_0$ and $z_0=i\sqrt{d_0}$ otherwise (where $\mQ(\sqrt{-d_0})=E$ and $d_0$ is a square free positive integer). In particular, any choice of $u$ has 
the form
$$
\textrm{(root of unity in $K$)}\cdot{ \Im(z_E)^{6h^2}|\eta^{24}(z_E)|^{h^2}(4\pi)^{6h^2}\over D^{-3h^2}\big[\prod_{j=1}^{D-1}\Gamma({j\over D})^{\legendre{-D}{j}}\big]^{{3w_Eh}}}.
$$
Secondly, note that by construction we have
have
$$
(\mFa,H/E)(u_c)=u_c\cdot\rho(\mFa^{-1})^{h}
$$
for any fractional ideal $\mFa$ of $E$. It was shown in \cite[chap. 2, p. 80]{Siegel-Adv} that the set of elements 
$\{\rho(\mFa)\}_{\mFa\in\Cl(E)\backslash\{0\}}$ is a multiplicatively 
independent set of units of $\CO_H$. Hence the orbit of the 
element $u_c\in H$ under $\Gal(H|E)$ has cardinality $\Gal(H|E)=\#\Cl(E)$, so 
$u_c$ is a class invariant of $E$, ie we have $H=E(u_c)$.  Numerical experiments suggest that the minimal polynomial of $u_c$ over $\mQ$ has in general small coefficients compared with those of 
$j(\CO_E)$ (which are famously large), so $u_c$ seems to be a "good" class invariant.

\end{document}